\documentclass[11pt, a4paper]{amsart}

\usepackage[a4paper, top=2.5cm, bottom=2.5cm, left=2.5cm, right=2.5cm]{geometry}

\usepackage{iftex}

\ifTUTeX
  \usepackage{fontspec}
  \usepackage[english]{babel}
  \babelprovide[import, onchar=ids fonts]{english}
  \babelfont{rm}{Noto Sans}
  \babelfont{sf}{Noto Sans}
\else
  \usepackage[utf8]{inputenc}
  \usepackage[T1]{fontenc}
  \usepackage{lmodern} 
  \usepackage[english]{babel}
\fi

\usepackage{amsmath}
\usepackage{amssymb}
\usepackage{amsthm}
\usepackage{mathtools} 
\usepackage{mathrsfs}  

\usepackage{xcolor}

\usepackage{enumitem}
\setlist[itemize]{label=-}
\usepackage{hyperref}
\hypersetup{
    colorlinks=true,
    linkcolor=blue,
    citecolor=red,
    urlcolor=blue
}
\usepackage{url}

\newcommand{\Z}{\mathbb{Z}}

\newcommand{\Primes}{\mathcal{P}}  
\newcommand{\A}{\mathcal{A}}         

\theoremstyle{plain}
\newtheorem{theorem}{Theorem}
\newtheorem{lemma}{Lemma}

\newtheorem{conjecture}{Conjecture}

\newtheorem*{theoremA}{Theorem A}

\theoremstyle{definition}

\title{\textbf{Almost-primes in Sun's $\lowercase{x^2+ny^2}$ Conjecture}}
\author{Songlin Han}
\address{Graduate School of Mathematics, Kyushu University, Fukuoka, Japan.}
\email{han.songlin.638@s.kyushu-u.ac.jp}

\author{Jinbo Yu}
\address{Graduate School of Mathematics,  Nagoya University, Nagoya, Japan.}
\email{jinbo.yu.e6@math.nagoya-u.ac.jp}

\subjclass[2020]{11N32, 
11N36
}
\keywords{linear sieve, weighted sieve, integer partitions, almost prime numbers}

\date{\today}

\let\origmaketitle\maketitle
\def\maketitle{
  \begingroup
  \def\uppercasenonmath##1{} 
  \let\MakeUppercase\relax 
  \origmaketitle
  \endgroup
}

\begin{document}

\begin{abstract}

\vspace{0.5em}
In 2015 Zhi-Wei Sun proposed the conjecture that any integer $n > 1$ admits a partition $n = x + y$ with integers $x, y >0$ such that $x + ny$ and $x^2 + ny^2$ are simultaneously prime. To approach this conjecture we use the method of weighted sieve as developed by Richert, Halberstam, and Diamond.
In this article, we first formalize the conjecture into a sieve problem.
We verify that the conditions required to use Richert's weighted sieve are satisfied and establish partial results with almost-prime solutions for sufficiently large $n$.
\end{abstract}

\maketitle

\section{Introduction}

In 2015, Zhi-Wei Sun introduced a series of conjectures \cite{Sun2017}, one of which asserts the following:

\begin{conjecture}[{\cite[Conjecture 2.21 (i)]{Sun2017}}]
\label{conj:A}
   Every integer $n > 1$ has a decomposition $n = x + y$ with integers $x,y > 0$, such that $x + ny$ and $x^2 + ny^2$ are simultaneously prime.
\end{conjecture}

In this paper, we approach Sun's problem using weighted sieve methods.
Sieve theory is one of the oldest methods in number theory.
It originated with the Sieve of Eratosthenes as an algorithm and was later formalized by Legendre.
Following significant improvements by Brun and Selberg, the theory entered a period of rapid development.
Although there are still limitations and unsolved problems within the theory, it remains a powerful tool for detecting almost-primes.
In particular, Richert proved that for any irreducible polynomial $F$ with degree $g = \deg F \ge 1$, there exist infinitely many integers $n$ such that $\Omega(F(n))\leq g+1$, where $\Omega(n)$ denotes the number of prime factors counted with multiplicity (see \cite[Theorem 6]{Richert1969}).

We apply Richert's method to study Sun's conjecture. Our result is as follows.

\begin{theorem}\label{main1}
    Every sufficiently large integer $n$ can be represented
    \begin{enumerate}[label=(\roman*),ref=(\roman*)]
        \item \label{thm11} as a sum of two positive integers $n = x + y$, such that \[\Omega(x + ny) \le 3,\]
        or
        \item \label{thm12} as a sum of two positive integers $n = x + y$, such that \[\Omega(x^2 + ny^2) \le 4.\]
    \end{enumerate}
\end{theorem}

\bigskip

Since $n = x+y$ implies $x = n-y$, \ref{thm11} and \ref{thm12} in Theorem~\ref{main1} are equivalent to  
 \begin{enumerate}[label=(\roman*),ref=(\roman*)]
    \item \label{thm11re} For any sufficiently large integer $n$  there is a positive integer $y$
    such that
    \begin{equation*}
       \Omega( n+(n-1)y)\leq3 .
    \end{equation*}
    \item \label{thm12re} For any sufficiently large integer $n$ there is a positive integer $y$
    such that
    \begin{equation*}
        \Omega((n+1)y^2 -2ny + n^2) \leq 4.
    \end{equation*}
    \end{enumerate}

If we focus only on part \ref{thm11re} of this problem, namely the existence of a prime of the form
\begin{align*}
    n+(n-1)y,\qquad 0<y<n,
\end{align*}
this may be viewed as a special case of Linnik's problem concerning upper bounds for the least prime in an arithmetic progression:
\begin{theoremA}[\cite{Linnik1944}]
    For integers $k,l >0$, with $l$ coprime to $k$, let $p(k,l)$ denote the smallest prime congruent to $l$ modulo $k$. Then for some fixed constant $L>0$, we have
    \begin{align*}
        \max_{(l,k)=1}p(k,l)\ll k^{L}.
    \end{align*}
\end{theoremA}
Recently, Xylouris \cite{Xylouris2011} showed that $L\leq5.18$, following the method of Heath-Brown. Chowla \cite{Chowla1934} observed that the Generalized Riemann Hypothesis implies $\max_{(l,k)=1}p(k,l)\ll_{\varepsilon} k^{2+\varepsilon}$ for any fixed $\varepsilon>0$ and conjectured that $\max_{(l,k)=1}p(k,l)\ll_{\varepsilon} k^{1+\varepsilon}$.
On the other hand, in 2016, Li--Pratt--Shakan \cite{LPS2017} established the following strong lower bounds
$$\max_{(l,k)=1}p(k,l)\gg\phi(k)(\log k)(\log\log k)\frac{\log\log\log\log k}{\log\log\log k},$$
with $\phi(k)$ denoting the Euler totient function,
indicating that for almost all moduli $k$, the first prime in the worst residue class is significantly larger than expected. These results indicate that prime-level statements in such problems are highly subtle, and naturally motivate the consideration of almost-prime variants.

Using a variation of Richert's method, known as the multi-dimensional weighted sieve (or the Diamond-Halberstam-Richert sieve, see \cite[Chapter 11]{DH2008}), we can prove that any sufficiently large integer $n$ can be represented as a sum of two positive integers $n=x+y$ such that, 
\begin{align*}
    \Omega\left((x+ny)(x^2+ny^2)\right) \le 11.
\end{align*}

\medskip

In this article, we use both Landau symbol and Vinogradov symbol for asymptotic analysis.

\section{Preliminaries}
\label{sec:pre}
In this section, we introduce some sieve notation and background information.

Let 
\begin{equation*}
    \A_n^i:=\{ F^i_n(y) : 1\le y < n\}, \quad (i =1,2)
\end{equation*}
and 
$F_n^{1}(y):=n+(n-1)y$, $F_n^2(y):=(n+1)y^2 -2ny + n^2$.
Let $X_i:= \# \A_n^i$ for $i=1,2$, 
thus $X_i = n-1$ for sufficiently large $n$. 

We use general notation in sieve theory. We denote by $\Primes$ the set of all primes and set
\begin{equation*}
    P(z):=\prod_{\substack{p\in\Primes\\p<z}}p
\end{equation*}
for $z\ge 2$.
In this article, the letter $p$ always denotes a prime, that is, $p\in\Primes$.
For integers $a,b>0$, $(a,b)$ denotes the greatest common divisor of $a$ and $b$.
For $\A \subset \Z_{>0}$ and $z\geq 2$ we define 
\begin{equation*}
    S(\A,z):=|\{a\in\A:(a,P(z))=1\}|
\end{equation*}
and for some square-free positive integer $q$ with $(q,P(z))=1$, we let
\begin{equation*}
    S_{q}(\A,z):=|\{a\in\A:a\equiv 0 \bmod q,(a,P(z))=1\}|.
\end{equation*}

For a prime $p$, we define the local density function $\rho(\A,p)$ as the number of solutions to the congruence equation $F_n(y) \equiv 0 \pmod p$ in $\A$. Here, $F_n$ denotes either $F^1_n$ or $F^2_n$.
Then $\rho(\A,d)$ can be defined by the Chinese remainder theorem for any square-free number $d$.
Clearly,
\begin{equation}\label{A1}
    0\le \rho(\A,p)\le \left(1-\frac{1}{A_1}\right)p
\end{equation}
for some constant $A_1\ge 1$. This follows from the observation that any polynomial of degree $g$ has at most $g-1$ solutions if it has no fixed prime factor, or in other words, the greatest common divisor of all coefficients of the polynomial is $1$.
By the definition of $\rho$, we also have 
\begin{equation*}
    \sum_{\substack{a\in\A\\a\equiv 0\bmod d}} 1
    =\sum_{\substack{1\le l\le d\\F_n(l)\equiv 0 \bmod d}}\, \sum_{\substack{m\le n-1\\m\equiv l\bmod d}}1
    =\frac{\rho(\A,d)}{d}(n-1)+r\rho(\A,d),
\end{equation*}
where $|r|< 1$.
For $w\ge 2$, Mertens' theorem gives us
\begin{equation*}
    \sum_{p<w}\frac{\log p}{p}=\log w +O(1).
\end{equation*}
On the other hand, Nagell's result \cite[Equation (4)]{Nagel1921} asserts that
\begin{equation*}
    \sum_{p<w}\frac{\rho(\A,p)}{p}\log p=\log w +O(1).
\end{equation*}
Putting these together, we obtain
\begin{equation}\label{A2}
    \Bigg| \sum_{z\le p < w} \frac{\rho(\A,p)-1}{p}\log p\Bigg|\le A_2
\end{equation}
for $w\ge z\ge 2$ and some constant $A_2\ge 1$.
For $z\ge 2$, we define the global density function by a finite product of the local density function,
\begin{equation*}
    G(z):= \prod_{p<z}\left(1-\frac{\rho(\A,p)}{p}\right).
\end{equation*}
Since neither $F_n^1$ nor $F_n^2$ has a fixed prime divisor, we have $\rho(\A,p) < p$ for all $p$, which ensures that $G(z) > 0$.

Now we consider the remainder term
\begin{equation*}
    \eta(\A,d)=\Bigg|\sum_{\substack{a\in\A\\a\equiv 0\bmod d}}1-\frac{\rho(\A,d)}{d}X \Bigg|
\end{equation*}
where $d$ is still a square-free number, $X=\# \A$ is the size of $\A$.
By the definition of $\eta(\A,d)$, it is clear that 
\begin{equation*}
    \eta(\A,d)\le \rho(\A,d)\le g^{\omega(d)}
\end{equation*}
where $g:=\deg F_n$ is $1$ if $F_n=F_n^{1}$ and $2$ if $F_n=F_n^{2}$,
and $\omega$ is the prime divisor function (without multiplicity).
Therefore, we see that 
\begin{equation*}
    \sum_{d\le X/\log^{3g+2} X}\mu^2(d)3^{\omega(d)}\eta(\A,d)
    \le \sum_{d\le X/\log^{3g+2} X}\mu^2(d)3^{\omega(d)}g^{\omega(d)}
    = \sum_{d\le X/\log^{3g+2} X}\mu^2(d)(3g)^{\omega(d)}.
\end{equation*}
Thus \cite[Lemma 3]{Richert1969} gives
\begin{equation}\label{A3}
    \sum_{d\le X/\log^{3g+2} X}\mu^2(d)(3g)^{\omega(d)}\le A_3 \frac{X}{\log^{3g+2}X}\log^{3g}X
\end{equation}
for some $A_3\ge 1$.
Also, for any fixed prime $p$, \cite[Th\'{e}or\`{e}me II]{Nagel1921} gives
\begin{equation}\label{A44}
    \sum_{\substack{a\in\A\\a\equiv 0 \bmod p^2}}1\le
    A_4\left(\frac{X\log X}{p^2}+1\right)
\end{equation}
for some $A_4\ge 1$.
Richert \cite{Richert1969} proved a weighted inequality under the conditions \eqref{A1}, \eqref{A2}, \eqref{A3} and \eqref{A44}. This inequality is the core idea of proving our main theorem. We present our proof in the next section.

The following two lemmas describe the asymptotic behavior of the density functions.
\begin{lemma}[{\cite[Lemma 1]{Richert1969}}]
    We have
    \begin{equation*}
        \sum_{z\le p<w}\frac{\rho(\A,p)}{p}=\log\frac{\log w}{\log z}+O\left(\frac{1}{\log z}\right)\quad \text{for }w\ge z,
    \end{equation*}
    and 
    \begin{equation*}
        \sum_{z\le p<w}\frac{\rho(\A,p)}{p}\log p=\log w-\log z+O(1)\quad \text{for }w\ge z.
    \end{equation*}
    \label{lem:A}
\end{lemma}

\begin{lemma}[{\cite[Lemma 2]{Richert1969}}]\label{lemB}
    We have
    \begin{equation*}
        \frac{1}{G(z)}\ll \log z,
    \end{equation*}
    \begin{equation*}
        G(z)=\prod_{p}\frac{1-\frac{\rho(\A,p)}{p}}{1-\frac{1}{p}}\frac{e^{-\gamma}}{\log z}\left(1+O\left(\frac{1}{\log z}\right)\right)
    \end{equation*}
    and
    \begin{equation*}
        \prod_{p}\frac{1-\frac{\rho(\A,p)}{p}}{1-\frac{1}{p}}\ge B_1
    \end{equation*}
    for some $B_1>0$.
    \label{lem:B}
\end{lemma}

Applying Selberg's sieve, we can show the following result.
See \cite{HR1974} or \cite{Richert1969} for the details.
\begin{lemma}\label{lemC_Selberg_sieve}
    If $\xi>B$ and if $z\ll \xi ^{A}$, we have
    \begin{equation*}
        \begin{split}
            S_q(\A,z)\le \frac{\rho(\A,q)}{q}XG(z)
                F\left(\frac{\log \xi^2}{\log z}\right)
            +\sum_{\substack{n\le \xi ^2\\n|P(z)}}3^{\omega(n)}\eta(\A,qn),\\
            S_q(\A,z)\ge \frac{\rho(\A,q)}{q}XG(z)
                f\left(\frac{\log \xi^2}{\log z}\right)
            -\sum_{\substack{n\le \xi ^2\\n|P(z)}}3^{\omega(n)}\eta(\A,qn).
        \end{split}
    \end{equation*}
    Furthermore, for $z\le X$
    \begin{equation*}
        \begin{split}
            S(\A,z)\le XG(z)F\left(\frac{\log X}{\log z}\right),\\
            S(\A,z)\ge XG(z)f\left(\frac{\log X}{\log z}\right).
        \end{split}
    \end{equation*}
    Here the functions $F$ and $f$ satisfy
    \begin{equation*}
        F(u)=\frac{2e^{\gamma}}{u}, \quad f(u)=0,\quad \text{ for }0<u\le 2,
    \end{equation*}
    and the differential equations
    \begin{equation*}
        (uF(u))'=f(u-1),\quad (uf(u))'=F(u-1),\quad\text{for }u\ge 2,
    \end{equation*}
    where $\gamma$ is the Euler constant.
\end{lemma}

Next, we introduce some important properties of $F$ and $f$.
For $2\leq u\leq4$, we can calculate $F$ explicitly by its definition and obtain
\begin{equation*}
    f(u)=\frac{2{e^{\gamma}}}{u}\log(u-1).
\end{equation*}
Also, $F(u)$ is monotonically decreasing toward $1$ and conversely $f(u)$ is monotonically increasing toward $1$.
Furthermore, we recall that 
\begin{equation*}
    F(u_1)-F(u_2)\ll\frac{u_2-u_1}{u_2^2},
    \quad
    f(u_1)-f(u_2)\ll\frac{u_2-u_1}{u_1^2}
\end{equation*}
for $0<u_1<u_2$.

\section{Proof of Theorem \ref{main1}}
We consider Richert's weighted sifting function.
Following \cite{Richert1969}, we let 
\begin{equation}
    W(\A,  u, \lambda) := \sideset{}{^{\prime}}\sum_{\substack{a\in\A\\(a, P(X^{1/4}))=1}} \left( 1 - \lambda \sum_{\substack{{X^{1/4} \le p < X^{1/u}}\\p|a}} \left( 1 - u \frac{\log p}{\log X} \right) \right),
    \label{def_W}
\end{equation}
where $\sideset{}{^{\prime}}\sum$ indicates that the summation is restricted to those $a$'s for which 
\begin{equation*}
    a\not\equiv 0\bmod p^2 \quad\text{for}\quad  X^{1/4}\le p<X^{1/u}.
\end{equation*}

Crucially, this definition provides a lower bound for the number of almost-primes.
To see this, we bound \eqref{def_W} from above as follows:
\begin{equation*}
    W(\A,  u, \lambda)\le \sideset{}{^{\prime}}\sum_{\substack{a\in\A\\(a, P(X^{1/4}))=1}}
    \left(1-\lambda\left(\Omega(a)-u\frac{\log a}{\log X}\right)\right).
\end{equation*}
Now, if we define 
\begin{equation*}
    \Lambda_r := r+1-\frac{\log\frac{4}{1+3^{-r}}}{\log 3},
\end{equation*}
then for $r\ge 2$ we have
\begin{equation*}
    r+1 -\frac{\log 4}{\log 3}\le \Lambda_r\le r+1 -\frac{\log 3.6}{\log 3}
\end{equation*}
since $\frac{4}{1+3^{-r}}$ is $4$ as $r$ tends to infinity and $3.6$ if $r=2$.
Furthermore, if we take
    $a\le X^{\Lambda_r-\delta}$
for all $a\in\A$ and some real $\delta\in \left(0,\frac{2}{3}\right]$,
then
\begin{equation}\label{Lambda}
    \frac{\log a }{\log X}\le \Lambda_r-\delta
\end{equation}
for sufficiently large $X$.
Thus we have
\begin{equation*}
    W(\A,  u, \lambda)\le 1-\lambda\left(\Omega(a)-u(\Lambda_r -\delta)\right).
\end{equation*}
Taking $u=1+3^{-r}$ and $\frac{1}{\lambda}=r+1-u(\Lambda_r -\delta)$, we easily see that 
\begin{equation}
    W(\A,  u, \lambda)\le \sum_{\substack{a\in\A \\ \Omega(a)\le r \\(a,P(X^{1/4}))=1}} 1\le \sum_{\substack{a\in\A \\ \Omega(a)\le r }} 1.
\end{equation}
\medskip

On the other hand, \cite{Richert1969} gives the following lower bound for $W(\A,  u, \lambda)$:
\begin{equation}\label{richert_theorem_1}
    W(\A,u,\lambda)\ge XG(X^{1/4})\left(\frac{e^{\gamma}}{2}\log 3-\lambda\int_{u}^{4}F\left(4\left(1-\frac{1}{t}\right)\right)\left(1-\frac{u}{t}\right)\frac{dt}{t}-\frac{b}{(\log X)^{1/15}}\right),
\end{equation}
where $\gamma$ is the Euler constant and $b$ is a constant depending only on $u$. We briefly outline the proof of this inequality below; for full details, we refer the reader to \cite[Theorem 1]{Richert1969}.

The strategy is to relate $W$ to the unweighted sifting function $S(\A, z)$ and then subtract error terms corresponding to elements with prime factors in some range $[z, y)$.
If we let $z=X^{1/4}$ and $y=X^{1/u}$, then by definition \eqref{def_W} we have
\begin{equation}
\label{eq:3term}
    W(\A,  u, \lambda)\ge S(\A,z)-\sum_{z\le p<y}\sum_{\substack{a\in\A\\a\equiv0\bmod p^2}}1-\lambda\sum_{z\le p<y}\left(1-\frac{\log p}{\log y}\right)S_{p}(\A,z).
\end{equation}

 The first term on the right-hand side is 
    \begin{equation}\label{a}
        S(\A,z)\ge XG(z)\left(f(4)-B\frac{\log\log X}{(\log X)^{1/14}}\right).
    \end{equation}
 By Lemma \ref{lemB}, we have $\frac{1}{G(z)}\ll \log z$, and thus
 applying Lemma \ref{lemC_Selberg_sieve}, we get the following bound for the second term in \eqref{eq:3term}  
    \begin{equation}\label{aa}
        \begin{split}
            \sum_{z\le p<y}\sum_{\substack{a\in\A\\a\equiv0\bmod p^2}}1&\le A_4\sum_{z\le p<y}\left(\frac{X\log X}{p^2}+1\right)\\
            &\le A_4\left(\frac{X\log X}{z}+y\right)\\
            &\le XG(z)\frac{b}{\log X}
        \end{split}
    \end{equation}
 if 
    \begin{equation}\label{A4}
        \sum_{\substack{a\in\A\\a\equiv 0\bmod p^2}}1\le A_4 \left(\frac{X\log X}{p^2}+1\right)
    \end{equation}
    for some $A_4\ge 1$.
    The estimation for the last term in \eqref{eq:3term} is more involved.
    First, by Selberg's sieve, we have that the last term is at most a multiple of $\lambda$ of 
    \begin{equation*}
        \sum_{z\le p <y}\left(1-\frac{\log p}{\log y}\right)\left(\frac{\rho(\A,p)}{p}XG(z)\left(F\left(\frac{\log(\xi^2/p)}{\log z}\right)+b\frac{\log\log X}{(\xi/p)^{1/14}}\right)+\sum_{\substack{n\le \xi^2/p\\n|P(z)}}3^{\omega(n)}\eta(\A,pn)\right),
    \end{equation*}
    where $\xi^2=\frac{X}{(\log X)^B}$.
    By the properties of $F$, we know this is 
    \begin{equation*}
        \le XG(z)\left(\sum_{z\le p<y}\left(1-\frac{\log p}{\log y}\right)\frac{\rho(\A,p)}{p}F\left(\frac{\log(X/p)}{\log z}\right)+b\frac{\log\log X}{(\log X)^{1/14}}\right)
        +\sum_{d\le \xi^2}\mu^2(d)3^{\omega(d)}\eta(\A,d),
    \end{equation*}
    because 
    \begin{equation*}
        \sum_{d\le X/\log ^B X}\mu^2(d)3^{\omega(d)}\eta(\A,d)\le A_3\frac{X}{\log ^{15/14}X}
    \end{equation*}
    for some constant $B>0$.
    Then by Lemma~\ref{lemB} and Abel's identity, we get the following upper-bound for the last sum in \eqref{eq:3term}
    \begin{equation}\label{aaa}
      \lambda \sum_{z\le p<y}\left(1-\frac{\log p}{\log y}\right)S_{p}(\A,z)  \le \lambda\int_{u}^{4}F\left(4\left(1-\frac{1}{t}\right)\right)\left(1-\frac{u}{t}\right)\frac{dt}{t}+O\left(\frac{b\log\log X}{\log X}\right),
    \end{equation}
    which follows from \cite[Equation (3.10), (3.11)]{Richert1969}.
Combining \eqref{a}, \eqref{aa} and \eqref{aaa}, and using the fact $f(4)=\frac{e^{\gamma}}{2}\log 3$ we obtain \eqref{richert_theorem_1}, which is a special case of \cite[Theorem 1]{Richert1969}.

Noting that
\begin{equation*}
    4\left(1-\frac{1}{t}\right)\le 3
\end{equation*}
for $t\le 4$, and recalling the definition of $F$, we find that the integral in \eqref{richert_theorem_1} is 
\begin{equation*}
    \begin{split}
        \int_{u}^{4}\left(\frac{2e^{\gamma}}{4\left(1-\frac{1}{t}\right)}\right)\left(1-\frac{u}{t}\right)\frac{dt}{t}
        &=\frac{e^{\gamma}}{2}\int_{u}^{4}\frac{t}{t-1}\cdot\frac{t-u}{t}\cdot\frac{dt}{t}\\
        &=\frac{e^{\gamma}}{2}\int_{u}^{4}\frac{t-u}{t(t-1)}dt.
    \end{split}
\end{equation*}
Using integration by parts and the fact $\frac{t-u}{t(t-1)}=\frac{u}{t}+\frac{1-u}{t-1}$, we see that the integral above is 
\begin{equation*}
    =\frac{e^{\gamma}}{2}\left(u\log \frac{4}{u}-(u-1)\log\frac{3}{u-1}\right) =\frac{e^{\gamma}}{2}D(u),
\end{equation*}
where $D(u):=u\log \frac{4}{u}-(u-1)\log\frac{3}{u-1}$, which does not depend on $X$.
Substituting this into \eqref{richert_theorem_1}, we conclude 
\begin{equation*}
    \sum_{\substack{a\in\A\\\Omega(a)\le r}}1\ge W(\A,  u, \lambda)
    \ge XG(X^{1/4})\left(\frac{e^{\gamma}}{2}
    D(u)
    -\frac{b}{(\log X)^{1/15}}\right).
\end{equation*}

\bigskip

Now we are ready to prove Theorem \ref{main1}.
We begin with the first part of Theorem \ref{main1}, which concerns the linear form in Conjecture \ref{conj:A}. Recall that we aim to find a integer $0<y < n$ such that $\Omega(x + ny) \leq 3$, where $x = n-y$, so our target polynomial is $ F_n^1(y) = n+(n-1)y$.
Thus the sequence to be sifted is
\begin{equation*}
    \A_n^1 = \{ F_n^1(y) : 1 \le y < n \}.
\end{equation*}
From our previous arguments, we have
\begin{equation*}
    a=n +(n-1)y <n+(n-1)^2 =X_1^2+X_1+1
\end{equation*}
for all $a\in\A_n^1$ in this case.

To apply the weighted sieve result we have derived, we need to check condition \eqref{Lambda}, that is, we need to show
\begin{equation*} \label{eq:linear_cond}
    \frac{\log a}{\log X_1} \le \Lambda_r - \delta
\end{equation*}
for all $a\in\A_n^1$.

By using the definition of $\Lambda_r$, we see that $r=3$ is admissible because $\Lambda_3\approx 2.7$, and
\begin{equation*}
    \max_{a\in\A_n^1}\frac{\log a}{\log X_1}=\frac{\log (X_1^2+X_1+1)}{\log X_1}<2\frac{\log (X_1+1)}{\log X_1}\le \Lambda _3-\delta
\end{equation*}
for any sufficiently large $n$.

Therefore, applying the lower bound inequality derived earlier in this section, we have
\begin{equation*}
    \sum_{\substack{a\in\A_n^1\\\Omega(a)\le 3}} 1
    \ge X_1G(X_1^{1/4})\left(C-\frac{b}{(\log X_1)^{\frac{1}{15}}}\right)
\end{equation*}
for some constant $C > 0$ which does not depend on $X_1$.
Since the last factor is positive for sufficiently large $X_1$, the lower bound is strictly positive.
This implies the existence of at least one element $a \in \A_n^1$ such that $\Omega(a) \le 3$.
This completes the proof for case $F_n^1$.

Similarly, calculating the optimal $r$ in the condition \eqref{Lambda} for the case $F_n^2(y)$, we get that $r=4$ is admissible.
In this case, 
\begin{equation*}
    \A_n^2 = \{ F_n^2(y) : 1 \le y < n \}
\end{equation*}
and 
\begin{equation*}
    \sum_{\substack{a\in\A_n^2\\\Omega(a)\le 4}} 1
    \ge X_2G(X_2^{1/4})\left(C-\frac{b}{(\log X_2)^{\frac{1}{15}}}\right)>0
\end{equation*}
for sufficiently large $n$.
This completes the proof.

\section*{Acknowledgement}
The authors would like to thank Professor Ade Irma Suriajaya for her supervision. We also express our gratitude to Professor Kohji Matsumoto and Professor Henrik Bachmann for their kind advice. The second author was financially supported by JST SPRING, Grant Number JPMJSP2125, and would like to take this opportunity to thank the "THERS Make New Standards Program for the Next Generation Researchers".

\nocite{*}
\bibliographystyle{amsplain}
\bibliography{ref2}

\end{document}